\documentclass[11pt]{article}
\usepackage[T1]{fontenc}
\usepackage[utf8]{inputenc}
\usepackage{lmodern}
\usepackage{amsmath,amssymb,amsthm,mathrsfs}
\usepackage{mathtools}
\usepackage{amscd}
\usepackage[margin=1in]{geometry}
\usepackage{microtype}
\usepackage{parskip}
\usepackage{needspace}
\usepackage{xcolor}
\usepackage{enumitem}
\usepackage[hidelinks,bookmarksnumbered=true]{hyperref}

\newtheorem{theorem}{Theorem}[section]
\newtheorem{lemma}[theorem]{Lemma}
\newtheorem{corollary}[theorem]{Corollary}
\newtheorem{proposition}[theorem]{Proposition}
\theoremstyle{definition}
\newtheorem*{notations}{Notations}
\newtheorem*{conventions}{Conventions}
\theoremstyle{plain}
\setlist{itemsep=2pt,topsep=4pt}
\hypersetup{
 pdftitle={Equilateral-triangle-free sets in finite planes},
 pdfauthor={Daewoong Cheong, Thang Pham, Le Quang-Ham}}
\title{Equilateral-triangle-free sets in finite planes}
\author{Daewoong Cheong \and Thang Pham \and Le Quang-Ham}
\date{}

\begin{document}
\maketitle
\begin{abstract}
For a prime $p>3$, let $F(p)$ be the largest cardinality of a subset of
$\mathbb F_p^2$ containing no equilateral triangle with nonzero squared
side length. Alex Iosevich asked for the exact value of $F(p)$.
By constructing such subsets of $\mathbb F_p^2$ from cap sets in
$\mathbb F_3^n$, we obtain the uniform lower bounds
\[
F(p)\gg
\begin{cases}
p^{2\gamma},   & p\equiv 11\pmod{12},\\
p^{1+\gamma},  & p\equiv 1\pmod{12}.
\end{cases}
\]
Here $\gamma=\log_3(2.2202)$, so the two exponents are
$1.452008666\ldots$ and $1.726004333\ldots$, respectively.
We also apply a quantitative corners theorem to show that $F(p)=o(p^2)$
as $p\to\infty$ through either congruence class.
\end{abstract}

\section{Introduction}
For a prime $p>3$, let $\mathbb F_p$
denote the field with $p$ elements, and let $F(p)$ be the largest cardinality of a subset of
$\mathbb F_p^2$ containing no equilateral triangle with nonzero squared
side length. For {
$w=(a,b)\in\mathbb F_p^2$}, write
$\|w\|^2=a^2+b^2$. The \emph{squared distance} between
$w=(a,b)$ and $w'=(a',b')$ is
\[
 \|w-w'\|^2=(a-a')^2+(b-b')^2,
\]
and an \emph{equilateral triangle} consists of three distinct points
$X,Y,Z$ satisfying
\begin{equation}\label{eq:ordinary-equilateral}
 \|X-Y\|^2=\|Y-Z\|^2=\|Z-X\|^2\ne0.
\end{equation}
The problem of determining $F(p)$ was raised by Alex Iosevich.

One of the motivations behind this question comes from the \emph{finite-field simplex problem}.

The finite-field simplex problem asks how large a subset of a
finite vector space must be to determine many, or all, congruence classes
of simplices.
Two configurations are \emph{congruent} if a translation followed by an
orthogonal linear map takes one to the other.
This problem has been studied by Fourier analysis, incidence estimates,
group actions, and spectral graph methods{
~\cite{bhipr,bip,vinh-graphs, vinh-norm}.}
Bennett, Iosevich, and Pakianathan~\cite{bip} used the Elekes--Sharir
method to show that, for every prime $p\equiv 3\pmod{4}$, any set
$E\subseteq\mathbb F_p^2$ with $|E|\gg p^{7/4}$ determines
$\gg p^3$ congruence classes of triangles, a positive proportion of
all such classes. Bennett, Hart, Iosevich, Pakianathan, and
Rudnev~\cite[Theorem~1.6]{bhipr} {showed that, for every prime $p>3$,
the same conclusion holds whenever $|E|\gg p^{8/5}$.} For $p\equiv3\pmod{4}$, Pham~\cite[Theorem~1.2]{pham-fixed-side}
proved a stronger conclusion under the same size condition, namely, $E$ determines $\gg p^2$ triangle
congruence classes with any prescribed nonzero squared side length.

{
A set may determine many triangle congruence classes and still miss a
particular one. Suppose that $t_1,t_2,t_3$ are the
squared side lengths of a noncollinear triangle $\Delta$ in $\mathbb F_p^2$,
and that all three are nonzero. How large must $E\subseteq\mathbb F_p^2$
be to contain a triangle with these three squared side lengths?}

{
This question is better understood in higher dimensions. View $\Delta$ in
$\mathbb F_p^d$ by appending zero coordinates.} For $d\ge 4$, results of
Vinh~\cite{vinh-graphs} and Iosevich and Parshall~\cite{iosevich-parshall}
show that every set $E\subseteq\mathbb F_p^d$ with
$|E|\ge C p^{(d+3)/2}$ contains a congruent copy of $\Delta$.
Parshall's simplex theorem~\cite{parshall-simplices} gives the same
conclusion for $d\ge 3$ under the hypothesis $|E|\ge C p^{2(d+1)/3}$.
This last exponent equals $2$ at $d=2$, although the theorem requires
$d\ge 3$. {
In the plane, we know of no $\varepsilon>0$ for which
$|E|\ge C p^{2-\varepsilon}$ suffices for every such choice of
$t_1,t_2,t_3$.}

{
Consider first the case $t_1=t_2=t_3=1$.
How large can $E$ be if it contains no unit equilateral triangle?
We can also allow the common squared side length to vary and ask how
large $E$ can be if it contains no equilateral triangle of any nonzero
squared side length. This is Iosevich's question about $F(p)$.}

Bennett, Iosevich, and Pakianathan~\cite[Theorem~1.4 and Corollary~4.2]{bip}
showed that an equilateral triangle with nonzero squared side length exists in
$\mathbb F_p^2$ if and only if $3$ is a square in $\mathbb F_p$.
For $p>3$, this is equivalent to $p\equiv1,11\pmod{12}$.
Thus, $F(p)=p^2$ when $p\equiv5,7\pmod{12}$.
{
We therefore focus on $p\equiv1,11\pmod{12}$.}

The restriction to nonzero squared side length matters when $-1$ is a
square. In that case, the plane is partitioned into $p$ parallel
\emph{isotropic lines}, on each of which all squared distances are zero.
If triples of distinct points with all three squared
distances zero were also forbidden, each isotropic line would contain
at most two points of the set. The set would then have at most $2p$
points.

\subsection{Main results}
Polansky~\cite[Sections~3--4]{polansky}, in his honors thesis, constructed twelve points modulo
seven in \emph{triangular-lattice coordinates} and combined copies by base-seven
expansions. His lattice construction gives
\[
 F(p)\gg p^{\log_7 12},\qquad \log_7 12=1.276989\ldots,
\]
for the problem with nonzero squared side length considered here.

Using cap sets in $\mathbb F_3^n$, we improve this lower bound in both
residue classes.\footnote{In his concluding remarks, Polansky wondered whether there exists a connection with cap sets and, if so, whether it could improve his lower bound.}

\begin{theorem}\label{thm:lower}
Set $\gamma=\log_3(2.2202)$. There is an absolute constant $c>0$ such that, for every prime
$p\equiv1,11\pmod{12}$,
\[
 F(p)\ge c
 \begin{cases}
 p^{2\gamma},&p\equiv11\pmod{12},\\
 p^{1+\gamma},&p\equiv1\pmod{12}.
 \end{cases}
\]
\end{theorem}

Recall that a \emph{cap set} in $\mathbb F_3^n$ is a set
$C\subseteq\mathbb F_3^n$ for which
\[
 \xi+\eta+\zeta=0,\qquad \xi,\eta,\zeta\in C
 \quad\Longrightarrow\quad \xi=\eta=\zeta.
\]
Equivalently, $C$ contains no nonconstant three-term arithmetic
progression. The following proposition gives a lower bound for $F(p)$
in terms of $n$ and $|C|$, which is the key tool in proving Theorem \ref{thm:lower}.

\begin{proposition}\label{prop:cap-to-plane}
Let $n\ge1$, and let $C\subseteq\mathbb F_3^n$ be a nonempty
cap set. For every prime $p\equiv1,11\pmod{12}$,
\[
 F(p)\ge 6^{-2\log_3|C|}
 \begin{cases}
 p^{2\log_3|C|/n},&p\equiv11\pmod{12},\\
 p^{1+\log_3|C|/n},&p\equiv1\pmod{12}.
 \end{cases}
\]
\end{proposition}

The exponents in Theorem~\ref{thm:lower} are $2\gamma=1.452008666\ldots$ and
$1+\gamma=1.726004333\ldots$, respectively.
{
Since $1+\gamma>8/5$, our construction gives, for every sufficiently
large prime $p\equiv1\pmod{12}$, a set that determines $\gg p^3$
triangle congruence classes by~\cite[Theorem~1.6]{bhipr} yet contains
no equilateral triangle of nonzero squared side length.}

For the upper bound, we use a linear change of coordinates that sends
each corner with nonzero step to an equilateral triangle with nonzero
squared side length. Recall that a corner in $\mathbb F_p^2$ is a {triple} $(s,t), (s+h,t), (s,t+h)$ with $h\ne 0.$
The quantitative corners theorem of Jaber, Liu,
Lovett, Ostuni, and Sawhney~\cite{jllos} then gives the following bound.
In particular, $F(p)=o(p^2)$ as $p\rightarrow \infty$ through primes $p\equiv 1,11 \pmod {12}$.

\begin{theorem}\label{thm:upper}
There is an absolute constant $c>0$ such that, for every prime $p>3$
with $p\equiv1,11\pmod{12}$,
\[
 F(p)\le p^2\exp\!\left[-c(\log p)^{1/600}\right].
\]
\end{theorem}

{
\begin{corollary}
Let $c>0$ be as in Theorem~\ref{thm:upper}, and let $p>3$ be a
prime with $p\equiv1,11\pmod{12}$. Every set $E\subseteq\mathbb F_p^2$
satisfying
\[
 |E|>p^2\exp\!\left[-c(\log p)^{1/600}\right]
\]
contains an equilateral triangle with nonzero squared side length.
\end{corollary}

In the continuous setting, Maga~\cite{maga} constructed
a compact set $E\subseteq\mathbb R^2$ with $\dim_H(E)=2$ containing
no equilateral triangle. For $d\ge4$, Iosevich and
Liu~\cite{iosevich-liu} proved that compact subsets of $\mathbb R^d$
satisfying a Frostman bound with exponent sufficiently close to $d$
contain an equilateral triangle.
}

{
\textbf{Notation.}} We write $f\ll g$ or $g\gg f$ when $f\le c_0g$ for a constant
$c_0>0$. {
The implied constants are independent of the field size. We indicate any
other dependence explicitly.} We use $\log$ for the natural logarithm.

\Needspace{6\baselineskip}
\section{Proof of the key proposition}\label{sec:prime-transfer}
To prove Proposition~\ref{prop:cap-to-plane}, we first
describe the coordinate change and linear form used to construct equilateral-triangle-free sets. We then
construct integer point sets from cap sets and transfer them to
$\mathbb F_p^2$.

\subsection{Finite-field geometry}
Let $p>3$ be a prime. Let $\mathcal S=\{e_1=(1,0),\,e_2=(0,1)\}$ be
the standard basis of $\mathbb F_p^2$, and write $(a,b)$ for the
coordinates of a point with respect to $\mathcal S$. These are the
coordinates used in Section~1, so the squared length of a point is
$\|(a,b)\|^2=a^2+b^2$, and the dot product is
$(a,b)\cdot(a',b')=aa'+bb'$.

For the rest of this section, we assume that $p\equiv1,11\pmod{12}$.
By quadratic reciprocity, this is exactly the condition that $3$ is a
square in $\mathbb F_p$, and we fix $r\in\mathbb F_p$ with $r^2=3$.
Let
\[
 \mathcal B=\left\{v_1=(1,0),\ v_2=\left(\frac12,\frac r2\right)\right\}
\]
be a second ordered basis of $\mathbb F_p^2$, written in the coordinates
$(a,b)$, and write $(s,t)$ for the coordinates of a point with respect to
$\mathcal B$. The map
\[
 \Phi:\mathbb F_p^2\longrightarrow\mathbb F_p^2,\qquad
 \Phi(s,t)=\left(s+\frac t2,\frac{rt}{2}\right),
\]
converts $(s,t)$-coordinates into $(a,b)$-coordinates:
\[
 (a,b)=\Phi(s,t)=sv_1+tv_2 .
\]
It is invertible, since its determinant is $r/2\ne0$.
To distinguish the two coordinate systems, we write
\[
 \mathbb F_p^2=\bigl(\mathbb F_p^2,(s,t)\bigr),\qquad
 (\mathbb F_p^2)'=\bigl(\mathbb F_p^2,(a,b)\bigr),
\]
so that $\Phi:\mathbb F_p^2\to(\mathbb F_p^2)'$, and $(\mathbb F_p^2)'$ is
the plane of Section~1 in which $F(p)$ is defined.

We recall the criterion for equilateral triangles from~\cite{bip} and
verify the identities needed below. Define
\begin{gather*}
 R:\mathbb F_p^2\longrightarrow\mathbb F_p^2,\qquad
 R(s,t)=(-t,s+t),\\
 Q:\mathbb F_p^2\longrightarrow\mathbb F_p,\qquad
 Q(s,t)=s^2+st+t^2,\\
 T:(\mathbb F_p^2)'\longrightarrow(\mathbb F_p^2)',\qquad
 T(a,b)=\left(\frac{a-rb}{2},\frac{ra+b}{2}\right).
\end{gather*}
Direct computations show that 
for every $(s,t)\in\mathbb F_p^2$,
\[
 \Phi(R(s,t))=T(\Phi(s,t)),\qquad
 \|\Phi(s,t)\|^2=Q(s,t),\qquad
 Q(R(s,t))=Q(s,t).
\]

The formulas for $R$ and $Q$ have integer coefficients. We use the
same formulas to define $R:\mathbb Z^2\to\mathbb Z^2$ and
$Q:\mathbb Z^2\to\mathbb Z$, and also over $\mathbb R$. Over $\mathbb R$, the form $Q$ is positive definite
because $Q(s,t)=(s+t/2)^2+3t^2/4$.

\begin{notations}
Throughout, for a prime $p$, we identify $\mathbb F_p=\{\bar{0}, \ldots, \overline{p-1}\}.$
 For $m\geq 1$ and the rings $\mathbb F=\mathbb F_p, \mathbb Z,$ the map $I_m:\mathbb F^m\rightarrow \mathbb F^m$ denotes the identity map.  The map $\pi: \mathbb Z^m\rightarrow \mathbb F_p^m$ denotes the natural projection; and  $\theta: \mathbb F_p^m\rightarrow \mathbb Z^m$ denotes  a set-theoretic section of $\pi$, i.e., $\pi\circ \theta = I_m$ on $\mathbb F_p^m$, such that $\mathrm{Im} (\theta)=\{0,\ldots, p-1\}^{m}$. For a subset $B$ of $\mathbb F_p^m,$ let $\widetilde{B}:=\theta (B),$ which is referred to as a lift of $B$ (along $\pi$). For $v\in \mathbb Z^m,$ we denote $\bar{v}=\pi(v).$
\end{notations}

\begin{conventions}
The letters $X,Y,Z$ denote points of $(\mathbb F_p^2)'$, 
and the letters $x,y,z$ denote points of $\mathbb F_p^2$. The letters $\xi,\eta,\zeta$ denote points of
$\mathbb F_3^n$, $(\mathbb F_3^2)^n$, or $\mathbb Z^{2n}$, including
integer points of $\mathbb Z^2$. Scalars in $\mathbb F_p$
obtained from the linear form $\ell$ of
Lemma~\ref{geom:prime-geometry}\textup{(iv)} are denoted $u$,
$u_1,u_2,u_3$.  

\end{conventions}

\begin{lemma}\label{geom:prime-geometry}
{
{Let} $p>3$ be a prime. Then the following hold.
\begin{enumerate}[label=\textup{(\roman*)},leftmargin=*]
\item An equilateral triangle with nonzero squared side length exists
in $(\mathbb F_p^2)'$ if and only if $p\equiv1,11\pmod{12}$.
\end{enumerate}
Suppose that $p\equiv1,11\pmod{12}$, and let $r$, $\Phi$, $R$, $Q$,
and $T$ be as above.
\begin{enumerate}[label=\textup{(\roman*)},leftmargin=*,resume]
\item Three points $X,Y,Z\in(\mathbb F_p^2)'$ form an equilateral
triangle with nonzero squared side length if and only if, after
possibly interchanging $X$ and $Y$,
\[
 Z-X=T(Y-X)\qquad\text{and}\qquad\|Y-X\|^2\ne0.
\]
\item For $x,y,x\in \mathbb F_p^2$, the points $\Phi(x),\Phi(y),\Phi(z)\in (\mathbb F_p^2)'$ form an equilateral triangle with
nonzero squared side length if and only if, after
possibly interchanging $x$ and $y$,
\begin{equation}\label{eq:lattice-equilateral}
 z-x=R(y-x)\qquad\text{and}\qquad Q(y-x)\ne0.
\end{equation}

\item Suppose additionally that $p\equiv1\pmod{12}$.
There exists $\lambda\in\mathbb F_p$ satisfying
$\lambda^2-\lambda+1=0$. The linear form
\[
 \ell:\mathbb F_p^2\longrightarrow\mathbb F_p,\qquad
 \ell(s,t)=s+\lambda t,
\]
satisfies
\begin{equation}\label{eq:projection-identities}
 \begin{aligned}
 \ell(R(s,t))&=\lambda\ell(s,t),\\
 Q(s,t)&=\ell(s,t)(s+\lambda^{-1}t)
 \end{aligned}
\end{equation}
for every $(s,t)\in\mathbb F_p^2$.
\end{enumerate}
}
\end{lemma}

Parts (ii) and (iii) characterize equilateral triangles with nonzero squared side length by linear equations. The map \(T\) preserves squared lengths and has determinant \(1\), so we call it a rotation; over \(\mathbb R\), taking \(r=\sqrt3\) gives rotation through \(60^\circ\). We call (2) the rotation equation, which is a main tool in our argument.

Part~\textup{(iv)} provides a linear form whose \emph{fibers}
$\ell^{-1}(u)$, with $u\in\mathbb F_p$, are parallel affine lines in
$\mathbb F_p^2$. If $x,y\in\mathbb F_p^2$ satisfy $\ell(x)=\ell(y)$,
then~\eqref{eq:projection-identities} gives
\[
 Q(y-x)=0,\qquad \|\Phi(y)-\Phi(x)\|^2=0.
\]
Thus, $\Phi$ sends these fibers to parallel isotropic lines in
$(\mathbb F_p^2)'$.  Note that the map $\ell$ is defined only in the case $p\equiv 1 \pmod{12}$. Using this map, we will obtain a better bound in Proposition \ref{prop:cap-to-plane} in this case

\begin{proof}
Define $J(a,b)=(-b,a)$. Suppose that $X,Y,Z$ form an equilateral
triangle with nonzero squared side length $\tau$. Put
$W=Y-X$ and $W'=Z-X$. Then
\[
 \|W\|^2=\|W'\|^2=\tau,\qquad W\cdot W'=\tau/2.
\]
Since $\det(W,JW)=\tau\ne0$, the vectors $W,JW$ form a basis.
Using $W\cdot JW=0$, we can write
\[
 W'=\frac12W+\frac{\rho}{2}JW
\]
for some $\rho\in\mathbb F_p$. Direct computation show that 
\[
\|W'\|^2=\frac{1+\rho^2}{4}\|W\|^2.
\]
Thus we have
\[
 \tau=\frac{1+\rho^2}{4}\tau,
\]
so $\rho^2=3$. By quadratic reciprocity, this implies
$p\equiv1,11\pmod{12}$. Conversely, if $r^2=3$, then
$(0,0),(2,0),(1,r)$ form an equilateral triangle of nonzero
squared side length. This proves (i).

For (ii), direct computation gives
\begin{equation} \label{3equations}
 T=\frac12I_2+\frac r2J,\qquad T^{-1}=I_2-T,\qquad
 \|Tw\|^2=|w|^2, \qquad \|(T-I_2)w\|^2=\|w\|^2.
\end{equation}

Since $\rho=\pm r$, the calculation in (i) shows that either
$Z-X=T(Y-X)$ or $Z-X=T^{-1}(Y-X)$. In the latter case,
$T^{-1}=I_2-T$ gives
\[
 Z-Y=T(X-Y).
\]
Thus, after possibly interchanging $X$ and $Y$, we obtain the
required equation. Conversely, if $Z-X=T(Y-X)$ and
$\|Y-X\|^2\ne0$, then using \eqref{3equations} with $w=Y-X$, we can show that
$X,Y,Z$ form an equilateral triangle with nonzero squared
side length. This proves (ii).

For (iii), apply (ii) to $\Phi(x),\Phi(y),\Phi(z)$.
Since $\Phi$ is a linear bijection satisfying
$\Phi R=T\Phi$ and $\|\Phi(w)\|^2=Q(w)$, the conditions in (ii)
are equivalent to
\[
 z-x=R(y-x),\qquad Q(y-x)\ne0,
\]
after possibly interchanging $x$ and $y$. This proves (iii).

Finally, suppose that $p\equiv1\pmod{12}$. Choose
$\iota\in\mathbb F_p$ with $\iota^2=-1$, and set
$\lambda=(1+\iota r)/2$. Then
\[
 \lambda(1-\lambda)=\frac{1+r^2}{4}=1,
\]
so $\lambda^2-\lambda+1=0$ and $\lambda^{-1}=1-\lambda$.
For $\ell(s,t)=s+\lambda t$, we obtain
\[
 \ell(R(s,t))=-t+\lambda(s+t)
 =\lambda(s+\lambda t)=\lambda\ell(s,t),
\]
and
\[
 (s+\lambda t)(s+\lambda^{-1}t)
 =s^2+st+t^2=Q(s,t).
\]
This proves (iv), completing the proof.
\end{proof}

\subsection{Encoding a cap set by integer points}
{
We construct integer point sets for which the rotation equation
forces the three points to coincide. The construction has three
steps, and we state each as a separate lemma.

For $n\ge1$, let $\mathcal R$ act on $(\mathbb F_3^2)^n$ and on
$\mathbb Z^{2n}$ by applying the formula for $R$ to each coordinate pair:
\[
 \mathcal R\bigl((s_j,t_j)_{j=0}^{n-1}\bigr)
 =\bigl(-t_j,s_j+t_j\bigr)_{j=0}^{n-1}.
\]
For $n=1$, this is the map $R$ itself. We say that a subset $S$ of
$(\mathbb F_3^2)^n$ or of $\mathbb Z^{2n}$ is \emph{rotation-free} if
\[
 \zeta=(I_{2n}-\mathcal R)\xi+\mathcal R \eta,\qquad \xi,\eta,\zeta\in S
 \quad\Longrightarrow\quad \xi=\eta=\zeta,
\]
with arithmetic in $\mathbb F_3$ or in $\mathbb Z$, respectively.
Lemma~\ref{lem:mod-three} produces a rotation-free subset of
$(\mathbb F_3^2)^n$ from a cap set. Lemma~\ref{lem:base-three} lifts
it to integer vectors by base-three expansions, and
Lemma~\ref{lem:planar-packing} packs the integer vectors into the
plane by expansions in a larger base. Proposition~\ref{lem:integer-packing}
combines the three steps.

\begin{lemma}\label{lem:mod-three}
Let $n\ge1$, and let $C\subseteq\mathbb F_3^n$ be a cap set. Define
\[
 D=\left\{((s_j,t_j))_{j=0}^{n-1}\in(\mathbb F_3^2)^n:
 (s_j-t_j)_{j=0}^{n-1}\in C,\ (s_j)_{j=0}^{n-1}\in C\right\}.
\]
Then $D$ is rotation-free and $|D|=|C|^2$.
\end{lemma}
}
\begin{proof}
We take a change of {coordinates} for which $R$ is triangular over $\mathbb F_3$.
 For
\[
 \xi=\bigl((s_j,t_j)\bigr)_{j=0}^{n-1}\in(\mathbb F_3^2)^n,
\]
define
\[
 \alpha_\xi=(s_j-t_j)_{j=0}^{n-1}\in\mathbb F_3^n,\qquad
 \beta_\xi=(s_j)_{j=0}^{n-1}\in\mathbb F_3^n.
\]
Define
\[
 \sigma:\mathbb F_3^{2n}=(\mathbb F_3^2)^n\longrightarrow\mathbb F_3^n\times\mathbb F_3^n=\mathbb F_3^{2n},
 \qquad \sigma(\xi)=(\alpha_\xi,\beta_\xi).
\]
Then $\sigma$ is a linear isomorphism with inverse
\[
 \sigma^{-1}(\alpha,\beta)=\bigl((\beta_j,\beta_j-\alpha_j)\bigr)_{j=0}^{n-1},
 \qquad \alpha=(\alpha_j)_{j=0}^{n-1},\ \beta=(\beta_j)_{j=0}^{n-1}
 \in\mathbb F_3^n.
\]
Identifying $(\mathbb F_3^2)^n=\mathbb F_3^{2n}$  and $ (\mathbb F_3^n)^2=\mathbb F_3^{2n}, $ $\sigma$ makes a change of coordinates.
Let $\mathcal R'=\sigma\circ\mathcal R\circ\sigma^{-1}$, so that the
diagram
\[
\begin{CD}
 (\mathbb F_3^2)^n @>\sigma>> (\mathbb F_3^n)^2\\
 @V\mathcal R VV @VV\mathcal R' V\\
 (\mathbb F_3^2)^n @>\sigma>> (\mathbb F_3^n)^2
\end{CD}
\]
commutes. A direct computation using $3=0$ in $\mathbb F_3$ shows
\[
 \mathcal R'(\alpha,\beta)=(-\alpha,\alpha-\beta),
 \qquad(\alpha,\beta)\in(\mathbb F_3^n)^2.
\]
The set $D$ is the inverse image of $C\times C$ under $\sigma$:
\[
 D=\sigma^{-1}(C\times C)
 =\bigl\{\xi\in(\mathbb F_3^2)^n:\alpha_\xi\in C,\ \beta_\xi\in C\bigr\}
 =\bigl\{((\beta_j,\beta_j-\alpha_j))_{j=0}^{n-1}:
 \alpha,\beta\in C\bigr\}.
\]
Hence, $|D|=|C|^2.$

If $\zeta=(I_{2n}-\mathcal R)\xi+\mathcal R \eta$ with $\xi,\eta,\zeta\in D$, then
applying $\sigma$ and the formula for $\mathcal R'$ gives
\[
 \alpha_\zeta=\alpha_\xi-(-\alpha_\xi)+(-\alpha_\eta)=2\alpha_\xi-\alpha_\eta,\qquad
 \beta_\zeta=\beta_\xi-(\alpha_\xi-\beta_\xi)+(\alpha_\eta-\beta_\eta)
 =2\beta_\xi-\alpha_\xi+\alpha_\eta-\beta_\eta.
\]
Since $2=-1$ in $\mathbb F_3$, rearranging gives
\[
 \alpha_\xi+\alpha_\eta+\alpha_\zeta=0,\qquad
 \beta_\xi+\beta_\eta+\beta_\zeta=\alpha_\eta-\alpha_\xi.
\]
By the definition of a cap set, we have $\alpha_\xi=\alpha_\eta=\alpha_\zeta$, and hence $\beta_\xi=\beta_\eta=\beta_\zeta$.
Since $\sigma$ is injective, $\xi=\eta=\zeta$. Thus{,} $D$ is rotation-free, as desired.
\end{proof}


Let $B$ be any subset of $\left(\mathbb F_3^2\right)^{{n}}$, so that $\widetilde{B}$ is a subset of $\left(\mathbb Z^2\right)^n.$
For each $k\geq 1,$ define 
$\Gamma: \left(\widetilde{B}\right)^{k}\rightarrow \mathbb Z^{2n}$ by
\[
\Gamma(\delta_0,\ldots,\delta_{k-1})=\sum_{i=0}^{k-1}3^{i} \delta_i.
\]  Let $E_k$ be the image of $\Gamma$.  Then, for each component ${\nu_j}$ of an element $\nu$ of $ E_k$, we have $0\leq \nu_j \leq 3^k-1.$ Thus{,} $E_k\subseteq \{0, \ldots, 3^k-1\}^{2n}$.
 
For $\xi \in E_k$,  $\delta_{\xi,j}$  denotes $\delta_j\in \widetilde{B}$ for $j=0,\ldots,k-1$, where $\xi=\Gamma(\delta_0,\ldots,  \delta_{k-1})$. Note that $\xi\equiv \delta_{\xi,0}\pmod3$ by the definition of $\xi.$
{
\begin{lemma}\label{lem:base-three}
Let $n\geq1$, {let} $B$ {be} a subset of $(\mathbb F_3^2)^n$ and {let} $E_k$ {be} as above for $k\geq 1.$ Then the following hold.
\begin{enumerate}[label=\textup{(\roman*)},leftmargin=*]
\item For each $k\geq 1,$
$\Gamma$ is injective, and so $|E_k|=|B|^k.$
\item 
If $B$ is rotation-free, then so is $E_k$ in $ \mathbb Z^{2n}$ for every $k\geq 1.$
\end{enumerate}
\end{lemma}
}
\begin{proof} (i): Injectivity of $\Gamma$ follows from uniqueness of base-three expansions. Since $\Gamma$ is injective, we have $|E_k|=|B|^k.$

(ii): We proceed by induction on $k.$ 
Suppose that $\xi,\eta,\zeta\in E_k$ satisfy
\begin{equation} \label{rotation equation-integer}
 \zeta=(I_{2n}-\mathcal R)\xi+\mathcal R \eta.
\end{equation}

Since $\xi\equiv \delta_{\xi,0}{,}\quad \eta\equiv \delta_{\eta,0}, \, \, \zeta\equiv \delta_{\zeta,0}  \pmod3$, the {equation} \eqref{rotation equation-integer} {reduces to} 
\begin{equation} \label{rotation equation-mod3}
 \delta_{\zeta,0}\equiv(I_{2n}-\mathcal R)\delta_{\xi,0}
 +\mathcal R\delta_{\eta,0}\pmod3.
\end{equation}
Let $\bar{\delta}_{\xi,0},\bar{\delta}_{\eta,0},\bar{\delta}_{\zeta,0} \in B$ be images under $\pi$ of the elements $\delta_{\xi,0}, \delta_{\eta,0}, \delta_{\zeta,0} \in \widetilde{B}$.
By \eqref{rotation equation-mod3}, we have, over $\mathbb F_3,$
\begin{equation} 
 \bar{\delta}_{\zeta,0}=(I_{2n}-\mathcal R)\bar{\delta}_{\xi,0}
 +\mathcal R\bar{\delta}_{\eta,0}.
\end{equation}
Since $B$ is rotation-free, we have
\[
 \bar{\delta}_{\xi,0}=\bar{\delta}_{\eta,0}=\bar{\delta}_{\zeta,0}.
\]

Applying $\theta$ to this equality, we obtain
\[
 \delta_{\xi,0}=\delta_{\eta,0}=\delta_{\zeta,0},
\]
{since $\theta(\pi(\delta))=\delta$ for every $\delta\in\widetilde B$.}
If $k=1,$ which is the base case for the induction, then $E_1=\widetilde{B}$ and
\[
 \delta_{\xi,0}=\xi, \qquad \delta_{\eta,0}=\eta, \qquad \delta_{\zeta,0}=\zeta.
\] Thus{,} $\xi=\eta=\zeta,$ and so $E_1$ is rotation-free.

Now {assume} $k\geq 2$.
Write $\delta_0:=\delta_{\xi,0}=\delta_{\eta,0}=\delta_{\zeta,0}$. Define $\xi', \eta',\zeta'$ by
\[
 \xi=\delta_0+3\xi',\qquad \eta=\delta_0+3\eta',\qquad \zeta=\delta_0+3\zeta'.
\]
Substitute these into \eqref{rotation equation-integer}{:}
\[
 \delta_0+3\zeta'=(I_{2n}-\mathcal R)(\delta_0+3\xi')+\mathcal R (\delta_0+{3}\eta').
\]
 Then since $(I_{2n}-\mathcal R)+\mathcal R=I_{2n}$, we obtain
\[
 \zeta'=(I_{2n}-\mathcal R)\xi'+\mathcal R \eta'.
\] Note that  $\xi',\eta',\zeta'$ are elements of $E_{k-1}$. By {the} induction hypothesis that $E_{k-1}$ is rotation-free, we obtain
$\xi'=\eta'=\zeta'$. Thus{,} $\xi=\eta=\zeta$, and so $E_k$ is rotation-free.
\end{proof}

Let $N\geq 1$ {be an integer}. Define $\Psi: \left(\mathbb Z^2\right)^n \rightarrow \mathbb Z^2$ by 
\[
 \Psi(\xi):=\left(\sum_{j=0}^{n-1}N^j s_j,
 \sum_{j=0}^{n-1}N^j t_j\right)=\sum_{j=0}^{n-1}N^j (s_j,t_j),
\]
where $\xi=((s_j,t_j))_{j=0}^{n-1}\in\left(\mathbb Z^2\right)^n.$

For $\xi=((s_j,t_j))_{j=0}^{n-1}$, let $\xi^{(j)}$ be the $j$-th coordinate pair, i.e., $\xi^{(j)}=(s_j,t_j).$

\begin{lemma}\label{lem:planar-packing}
Let $n\geq 1$ and fix an integer $M\ge0$. Let $G$ be a subset of ${\mathbb Z^{2n}}$ with
$G\subseteq\{0,\ldots,M\}^{2n}.$
\begin{enumerate}[label=\textup{(\roman*)},leftmargin=*]
\item If $N\ge\max\{2,2M+1\}$, then $\Psi$ is injective on $G$, and 
\[
 \Psi(G)\subseteq\left\{0,\ldots,M\,\frac{N^n-1}{N-1}\right\}^2,
\]
\item If {$N\ge\max\{2,2M+1\}$ and} $G$ is rotation-free, then so is $\Psi(G)$  in $\mathbb Z^2$.
\end{enumerate}
\end{lemma}

\begin{proof}
(i): Since $0\le s_j,t_j\le M<N$ for $\xi\in G$, injectivity of $\Psi$ on $G$ follows from uniqueness of base-$N$
expansions. Both coordinates of
$\Psi(\xi)$ lie between $0$ and
\[
 M\sum_{j=0}^{n-1}N^j=M\,\frac{N^n-1}{N-1}.
\]

(ii): To prove (ii), 
assume that for $\xi,\eta,\zeta\in G$
\begin{equation}\label{rotation-ZZ}
\Psi(\zeta)=(I_2-R)\Psi(\xi)+R\Psi(\eta).
\end{equation} Let us show that $\Psi(\xi)=\Psi(\eta)=\Psi(\zeta).$
For each $j=0,\ldots{,}n-1,$ let $e^{(j)}\in \mathbb Z^2$ be the $j$-th error term  defined by
\[
 e^{(j)}=\zeta^{(j)}-(I_2-R)\xi^{(j)}-R\eta^{(j)}.
\]
Using the fact that $\Psi\circ \mathcal R= R\circ \Psi$, and $R, \mathcal R$ are linear, from \eqref{rotation-ZZ} we can derive  
 \begin{equation} \label{equation for e}
 \sum_{j=0}^{n-1}N^j e^{(j)}=0.
 \end{equation}
By rewriting the {equality} \eqref{equation for e} as \[e^{(0)}=(e^{(0)}_1,e^{(0)}_2)=-N\sum_{j=1}^{n-1}N^{j-1}e^{(j)},\] we can easily see that each component $e^{(0)}_i$ of $e^{(0)}$ ($i=1,2$) is {divisible} by $N$ in $\mathbb Z.$

On the other hand, each component  ${e^{(j)}_i}$ of $e^{(j)}$  is, by definition, a sum of two entries minus two entries for $i=1,2$,
and each entry lies between $0$ and $M$, so for every $j=0,\ldots,n-1,$
\[
 |e^{(j)}_1|,|e^{(j)}_2|\le2M\le N-1.
\]

For $j=0,$ since  $e^{(0)}_i$ is a multiple of $N$ and $|e^{(0)}_i|\leq N-1$ for $i=1,2$ we have $e^{(0)}_i=0$, so $e^{(0)}=(0,0).$ 

Removing this zero term from \eqref{equation for e} gives $\sum_{j=1}^{n-1}N^j e^{(j)}=0$. {After dividing by $N$, the same argument gives $e^{(1)}=(0,0)$ if $n\ge2$.} Continuing this process, we obtain 
\[e^{(0)}=e^{(1)}=\cdots =e^{(n-1)}=(0,0).\]
Thus{,} $\xi,\eta,\zeta \in G$ {satisfy} 
\[
 \zeta=(I_{2n}-\mathcal R)\xi+\mathcal R \eta.
\]
Since $G$ is rotation-free, we have $\xi=\eta=\zeta.$ Therefore
$\Psi(\xi)=\Psi(\eta)=\Psi(\zeta).$ Thus{,} the subset $\Psi(G)\subset \mathbb Z^2$ is rotation-free.

\end{proof}

We now combine the three lemmas.
\begin{proposition}\label{lem:integer-packing}
Let $n\geq 1$ and $C\subseteq\mathbb F_3^n$ {be} a nonempty cap set.
For every integer $k\ge0$, define $N_k=2\cdot3^k-1$.
There is an explicit set
\[
 S_k\subseteq\left\{0,\ldots,\frac{N_k^n-1}{2}\right\}^2,
 \qquad |S_k|=|C|^{2k},
\]
such that
\[
 \zeta=(I_2-R)\xi+R\eta,\qquad \xi,\eta,\zeta\in S_k
 \quad\Longrightarrow\quad \xi=\eta=\zeta.
\]
\end{proposition}
\begin{proof}

For $k=0$, take $S_0=\{(0,0)\}$. Assume henceforth that $k\ge1$.
Let $D$ be the rotation-free set of Lemma~\ref{lem:mod-three}, so that
$|D|=|C|^2$. Next, we apply Lemma \ref{lem:base-three} 
 to  $B=D$. Then we obtain the rotation-free subset  $E_k\subseteq\{0,\ldots,3^k-1\}^{2n}$ such that $|E_k|=|D|^k=|C|^{2k}$.  Lastly,  apply Lemma~\ref{lem:planar-packing} with $M=3^k-1$ and
$N=N_k=2M+1$, and set $S_k=\Psi(E_k)$. Then $|S_k|=|E_k|=|C|^{2k}$,
the set $S_k$ is rotation-free in $\mathbb Z^2$, and its coordinates
lie between $0$ and
\[
 M\,\frac{N_k^n-1}{N_k-1}=\frac{N_k^n-1}{2},
\]
since $N_k-1=2M$. This proves the proposition.
\end{proof}

\subsection{From integer points to the finite plane}
We now use the integer sets from
Proposition~\ref{lem:integer-packing} to construct equilateral-triangle-free
subsets of $\mathbb F_p^2$.

\begin{proof}[Proof of Proposition~\ref{prop:cap-to-plane}]
We use direct reduction modulo $p$ when $p\equiv11\pmod{12}$.
When $p\equiv1\pmod{12}$, we first construct a scalar set and then
use parallel isotropic lines to gain a factor of $p$.
We finish by choosing the largest admissible scale in each case.

Let $S_k$ be the set from Proposition~\ref{lem:integer-packing}, with
$N_k=2\cdot3^k-1$. Thus, $|S_k|=|C|^{2k}$, and
\[
 \zeta=(I_2-R)\xi+R\eta\quad\text{in }\mathbb Z^2,
 \qquad \xi,\eta,\zeta\in S_k
 \quad\Longrightarrow\quad \xi=\eta=\zeta.
\]
Let $\pi:\mathbb Z^2\to\mathbb F_p^2$ be reduction modulo $p$,
and write $\bar\xi=\pi(\xi)$ for $\xi\in\mathbb Z^2$.
Since $R$ has integer coefficients, $\pi(R\xi)=R\bar\xi$.
We will obtain equations modulo $p$ and use the coordinate bounds
to recover the exact integer equation above.
For $\xi,\eta,\zeta\in S_k$, define the integer vector
$e=\zeta-(I_2-R)\xi-R\eta$. Then
\par
\begin{align*}
 e_1&=\zeta_1-\xi_1-\xi_2+\eta_2,\\
 e_2&=\zeta_2+\xi_1-\eta_1-\eta_2,\\
 e_1+e_2&=\zeta_1+\zeta_2-\xi_2-\eta_1.
\end{align*}
Note that each of $e_1, e_2, e_1+e_2$ is expressed as a sum of two entries minus two entries.
Since every entry lies between $0$ and $(N_k^n-1)/2$, we obtain
\begin{equation}\label{eq:error-bound}
 |e_1|,\ |e_2|,\ |e_1+e_2|\le N_k^n-1.
\end{equation}

\emph{The case $p\equiv11\pmod{12}$.}
Choose $k$ with $N_k^n\le p$.
Define
\[
 E=\Phi(\pi(S_k))\subseteq(\mathbb F_p^2)'.
\]
Reduction modulo $p$ is injective on $S_k$, since its coordinates lie
between $0$ and $(p-1)/2$. As $\Phi$ is also invertible,
$|E|=|S_k|=|C|^{2k}$.

Suppose that $X,Y,Z\in E$ form an equilateral triangle with nonzero
squared side length. Choose the unique $\xi,\eta,\zeta\in S_k$ such that
\[
 X=\Phi(\bar \xi),\qquad Y=\Phi(\bar \eta),\qquad Z=\Phi(\bar \zeta).
\]
By~\eqref{eq:lattice-equilateral}, after possibly
interchanging $X,Y$ and the corresponding $\xi,\eta$, we have
\[
 \bar \zeta=(I_2-R)\bar \xi+R\bar \eta\qquad\text{in }\mathbb F_p^2.
\]
Thus, the integer error $e=\zeta-(I_2-R)\xi-R\eta$ satisfies
$\pi(e)=0$, so $p\mid e_1,e_2$. On the other hand,
\eqref{eq:error-bound} gives
\par
\[
 |e_1|,|e_2|\le N_k^n-1<p.
\]
Hence, $e=0$, and we have the exact equality
\[
 \zeta=(I_2-R)\xi+R\eta\qquad\text{in }\mathbb Z^2.
\]
Proposition~\ref{lem:integer-packing} now gives $\xi=\eta=\zeta$, so $X=Y=Z$,
a contradiction. Therefore,
\par
\begin{equation}\label{eq:direct-size}
 F(p)\ge|C|^{2k}\qquad\text{whenever }N_k^n\le p.
\end{equation}

\emph{The case $p\equiv1\pmod{12}$.} The argument made so far in the case $p\equiv11\pmod{12}$ applies verbatim when $p\equiv1\pmod{12}$. However, the case $p\equiv1\pmod{12}$ is more favorable since there are isotropic lines in ${(\mathbb F_p^2)'}$. Using this advantage we shall find a better bound in this case.

Choose $k$ with $N_k^n\le\sqrt p$.
Let $\lambda$ and $\ell$ be as in
Lemma~\ref{geom:prime-geometry}\textup{(iv)}. We first show that
\begin{equation}\label{eq:scalar-lift}
 \ell(\bar \zeta)=(1-\lambda)\ell(\bar \xi)+\lambda\ell(\bar \eta),
 \qquad \xi,\eta,\zeta\in S_k
 \quad\Longrightarrow\quad \xi=\eta=\zeta.
\end{equation}
Suppose that the scalar equation holds, and set
$e=\zeta-(I_2-R)\xi-R\eta\in\mathbb Z^2$.
The first identity in~\eqref{eq:projection-identities} gives
\[
 \ell(\bar e)=\ell(\bar \zeta)-(1-\lambda)\ell(\bar \xi)
 -\lambda\ell(\bar \eta)=0.
\]
The second identity then yields $Q(\bar e)=0$ in $\mathbb F_p$.
Equivalently, the integer $Q(e)$ is divisible by $p$.
We will show that $0\le Q(e)<p$, which forces $Q(e)=0$ over the integers.
\par
For real $s,t$,
\begin{equation}\label{eq:quadratic-bound}
 Q(s,t)\le\max\{s^2,t^2,(s+t)^2\}.
\end{equation}
Indeed, if $st\ge0$, then $Q(s,t)=(s+t)^2-st\le(s+t)^2$.
If $st<0$ and $|s|\ge|t|$, then $st+t^2\le0$, so $Q(s,t)\le s^2$;
the other case is symmetric.
Combining~\eqref{eq:error-bound} and~\eqref{eq:quadratic-bound}, we obtain
\[
 0\le Q(e)\le(N_k^n-1)^2<p.
\]
Since $p\mid Q(e)$, we have $Q(e)=0$ in $\mathbb Z$.
Positive definiteness of $Q$ over the integers gives $e=0$.
Thus, $\zeta=(I_2-R)\xi+R\eta$ in $\mathbb Z^2$, and
Proposition~\ref{lem:integer-packing} yields $\xi=\eta=\zeta$.
This proves~\eqref{eq:scalar-lift}.

The map $\xi\mapsto\ell(\bar \xi)$ is also injective on $S_k$.
Indeed, if $\ell(\bar \xi)=\ell(\bar \eta)$, then
\[
 \ell(\bar \eta)=(1-\lambda)\ell(\bar \xi)+\lambda\ell(\bar \xi).
\]
Then \eqref{eq:scalar-lift} applied to $(\xi,\xi,\eta)$ gives $\xi=\eta$. {This proves the injectivity.}
Consequently, the set
\[
 U=\{\ell(\bar \xi):\xi\in S_k\}\subseteq\mathbb F_p
\]
has cardinality $|U|=|S_k|=|C|^{2k}$. Moreover, we claim that for $u_1,u_2,u_3\in U $
\begin{equation}\label{eq:scalar-free}
 u_3=(1-\lambda)u_1+\lambda u_2 \qquad u_1,u_2,u_3 \in U
 \quad\Longrightarrow\quad u_1=u_2=u_3.
 \end{equation}
Indeed, choosing representatives $\xi,\eta,\zeta \in
S_k$ so that $\ell(\bar{\xi})=u_1, \ell(\bar{\eta})=u_2, \ell(\bar{\zeta})=u_3,$ and again applying~\eqref{eq:scalar-lift}, we obtain $u_1=u_2=u_3.$
Define
\[
 A=\ell^{-1}(U)=\{x\in\mathbb F_p^2:\ell(x)\in U\},\qquad E=\Phi(A)\subseteq(\mathbb F_p^2)'.
\]
For each $u\in U$, the fiber $\ell^{-1}(u)$ is the line
$\{(s,t)\in\mathbb F_p^2:s+\lambda t=u\}$, which contains exactly $p$ points.
The fibers are disjoint, and $\Phi$ is invertible, so
$|E|=|A|=p|U|={p|S_k|}=p|C|^{2k}$.

Suppose that $E\subseteq {(\mathbb F_p^2)'}$ contains an equilateral triangle $X,Y,Z$ with nonzero squared
side length. Assume that  $\Phi(x)=X, \Phi(y)=Y, \Phi(z)=Z$ for $x,y,z\in  A\subseteq \mathbb F_p^2.$ ( Note that $x,y,z$ are already points of $\mathbb F_p^2$, not elements of $S_k \subseteq \mathbb Z^2.$) 
Reordering $X,Y,Z$, if necessary, by~\eqref{eq:lattice-equilateral} we obtain
\[
 z=(I_2-R)x+Ry,\qquad Q(y-x)\ne0.
\]
Applying $\ell$ gives
\[
 \ell(z)=(1-\lambda)\ell(x)+\lambda\ell(y).
\]
Since $\ell(x),\ell(y),\ell(z)\in U$, by~\eqref{eq:scalar-free}
we obtain 
\[\ell(x)=\ell(y)=\ell(z).
\]In particular,
$\ell(y-x)=0$, so~\eqref{eq:projection-identities} gives $Q(y-x)=0$,
a contradiction. Therefore,
\par
\begin{equation}\label{eq:split-size}
 F(p)\ge p|C|^{2k}\qquad\text{whenever }N_k^n\le\sqrt p.
\end{equation}

Let $L=p$ in~\eqref{eq:direct-size} and $L=\sqrt p$ in
\eqref{eq:split-size}. Choose the largest $k\ge0$ with $N_k^n\le L$.
Such a $k$ exists because $N_0=1$. Since $N_k=2 \cdot 3^k-1,$ the maximal choice of $k$ gives
\[
 L<(2\cdot3^{k+1}-1)^n<6^n3^{kn},
 \qquad 3^k>L^{1/n}/6.
\]
If $|C|>1$, raising the latter inequality to the power $2\log_3|C|$
gives
\begin{equation}\label{inequality}
 |C|^{2k}\ge6^{-2\log_3|C|}L^{2\log_3|C|/n}.
\end{equation}
If $|C|=1$, then both sides of \eqref{inequality} equal $1$, so it trivially holds.
Substituting \eqref{inequality} into~\eqref{eq:direct-size} and~\eqref{eq:split-size}, we obtain both bounds in the proposition.
\end{proof}

\Needspace{7\baselineskip}
\section{Proofs of Theorems~\ref{thm:lower} and \ref{thm:upper}}\label{sec:main-proofs}
We combine Proposition~\ref{prop:cap-to-plane} with the cap set
construction of Romera-Paredes et al.~\cite{funsearch} to obtain the
exponents in Theorem~\ref{thm:lower}.

\begin{proof}[Proof of Theorem~\ref{thm:lower}]
The construction of Romera-Paredes et al.~\cite[p.~472 and Supplementary
Information, Figure~C.9]{funsearch} gives an integer $n\ge1$ and
a cap set $C\subseteq\mathbb F_3^n$ with
$|C|\ge 2.2202^n=3^{n\gamma}$.
Fix such $n$ and $C$. Since $\log_3|C|/n\ge\gamma$,
Proposition~\ref{prop:cap-to-plane} gives the required bounds with
$c=6^{-2\log_3|C|}>0$.
Since $n$ and $C$ are fixed independently of $p$, the constant $c$
is absolute. This proves Theorem~\ref{thm:lower}.
\end{proof}

\Needspace{7\baselineskip}
We use $\Phi$ to map each corner to an equilateral triangle with
nonzero squared side length. The quantitative corners theorem then
bounds the size of an equilateral-triangle-free set.

\begin{proof}[{Proof of Theorem~\ref{thm:upper}}]
Let $E\subseteq(\mathbb F_p^2)'$ contain no equilateral triangle with
nonzero squared side length. Let $r$ and $\Phi$ be as in
Section~\ref{sec:prime-transfer}, so that
$\Phi^{-1}(E)\subseteq\mathbb F_p^2$ is the same set in the
coordinates $(s,t)$.
We claim that $\Phi^{-1}(E)$ doesn't contain any corners.
Suppose, {to the contrary}, that $\Phi^{-1}(E)$ contains a corner
\[
 (s,t),\qquad(s+h,t),\qquad(s,t+h),\qquad h\ne0.
\]
By linearity and the definition of $\Phi$, the differences between
these image points are
\begin{align*}
 \Phi(s+h,t)-\Phi(s,t)&=\Phi(h,0)=(h,0),\\
 \Phi(s,t+h)-\Phi(s,t)&=\Phi(0,h)
 =\left(\frac h2,\frac{rh}{2}\right),\\
 \Phi(s,t+h)-\Phi(s+h,t)&=\Phi(-h,h)
 =\left(-\frac h2,\frac{rh}{2}\right).
\end{align*}
The first vector has squared length $h^2$. Since $r^2=3$, the other
two vectors have squared lengths
\[
 \left(\frac h2\right)^2+\left(\frac{rh}{2}\right)^2
 =\left(-\frac h2\right)^2+\left(\frac{rh}{2}\right)^2
 =\frac{1+r^2}{4}h^2=h^2.
\]

Since $h\ne0$, the image points $ \Phi(s,t), \Phi(s+h,t), {\Phi(s,t+h)}$ in $ (\mathbb F_p^2)'$ form an equilateral triangle
with nonzero squared side length, contrary to the choice of $E$.
Thus, $\Phi^{-1}(E)$ doesn't contain any corner.
Applying the corners theorem of Jaber, Liu, Lovett, Ostuni, and
Sawhney~\cite[Theorem~1.1, full version]{jllos} to the additive group
$\mathbb F_p$
gives
\[
 |E|=|\Phi^{-1}(E)|\le p^2\exp\!\left[-c(\log p)^{1/600}\right].
\]
This proves Theorem~\ref{thm:upper}.
\end{proof}


\begin{thebibliography}{99}
\raggedright
\small
\setlength{\itemsep}{2pt}

\bibitem{bhipr}
M.~Bennett, D.~Hart, A.~Iosevich, J.~Pakianathan, and M.~Rudnev,
\emph{Group actions and geometric combinatorics in $\mathbb F_q^d$},
Forum Mathematicum, \textbf{29}(1) (2017), 91--110.

\bibitem{bip}
M.~Bennett, A.~Iosevich, and J.~Pakianathan,
\emph{Three-point configurations determined by subsets of
$\mathbb F_q^2$ via the Elekes--Sharir paradigm},
Combinatorica, \textbf{34}(6) (2014), 689--706.


\bibitem{iosevich-liu}
A.~Iosevich and B.~Liu,
\emph{Equilateral triangles in subsets of $\mathbb R^d$
of large Hausdorff dimension},
Israel Journal of Mathematics, \textbf{231}(1) (2019), 123--137.

\bibitem{iosevich-parshall}
{
A.~Iosevich and H.~Parshall,
\emph{Embedding distance graphs in finite field vector spaces},
Journal of the Korean Mathematical Society, \textbf{56}(6) (2019), 1515--1528.}

\bibitem{jllos}
M.~Jaber, Y.~P.~Liu, S.~Lovett, A.~Ostuni, and M.~Sawhney,
\emph{Quasipolynomial bounds for the corners theorem},
in Proceedings of the 2025 IEEE 66th Annual Symposium
on Foundations of Computer Science,
IEEE, 2025, 2027--2036; full version, arXiv:2504.07006.

\bibitem{maga}
P.~Maga,
\emph{Full dimensional sets without given patterns},
Real Analysis Exchange, \textbf{36}(1) (2010/11), 79--90.


\bibitem{parshall-simplices}
{
H.~Parshall,
\emph{Simplices over finite fields},
Proceedings of the American Mathematical Society, \textbf{145}(6) (2017),
2323--2334.}



\bibitem{pham-fixed-side}
T.~Pham,
\emph{Triangles with one fixed side-length, a Furstenberg-type problem,
and incidences in finite vector spaces},
Forum Mathematicum, \textbf{37}(2) (2025), 359--372.

\bibitem{polansky}
Z.~Polansky,
\emph{Equilateral triangles in vector spaces over finite fields},
honors thesis, University of Rochester, 2019.
\url{https://www.sas.rochester.edu/mth/undergraduate/honorspaperspdfs/zachpolansky19.pdf}.

\bibitem{funsearch}
B.~Romera-Paredes, M.~Barekatain, A.~Novikov, M.~Balog,
M.~P.~Kumar, E.~Dupont, F.~J.~R.~Ruiz, J.~S.~Ellenberg,
P.~Wang, O.~Fawzi, P.~Kohli, and A.~Fawzi,
\emph{Mathematical discoveries from program search with large language models},
Nature, \textbf{625}(7995) (2024), 468--475.

\bibitem{vinh-graphs}
L.~A.~Vinh,
\emph{On kaleidoscopic pseudo-randomness of finite Euclidean graphs},
Discussiones Mathematicae Graph Theory,
\textbf{32}(2) (2012), 279--287.



\bibitem{vinh-norm}
L.~A.~Vinh,
\emph{The solvability of norm, bilinear and quadratic equations over finite
fields via spectra of graphs},
Forum Mathematicum, \textbf{26}(1) (2014), 141--175.

\end{thebibliography}
\end{document}